\documentclass[3p,12pt]{elsarticle}
\usepackage[english]{babel}
\usepackage{amssymb}
\usepackage{amsthm}
\usepackage{amsmath}
\usepackage{empheq}
\usepackage{graphicx}
\usepackage{subcaption}
\usepackage{multicol}
\usepackage[utf8]{inputenc}

\usepackage{xcolor}
\newcommand{\correctcolor}{black} 
\newcommand{\fix}[1]{\textcolor{\correctcolor}{#1}}
\newcommand{\ii}{\mathrm{i}}

\begin{document}
	\begin{frontmatter}
		\title{Preconditioning for time-harmonic Maxwell's equations using the Laguerre transform}
		\author{Andrew V. Terekhov}
		\ead{andrew.terekhov@mail.ru}
		\address{Novosibirsk State Technical University, 630073, Novosibirsk, Russia}
		\begin{abstract}
			A method of numerically solving Maxwell’s equations is considered for modeling harmonic electromagnetic fields. The vector finite element  method makes it possible to obtain a physically consistent discretization of the differential equations. However, solving large systems of linear algebraic equations with indefinite ill-conditioned matrices is a challenge. The large size and fill-in of the system matrix limit the applicability of direct Gaussian elimination, since direct solvers require substantial memory and computational effort. The proposed iterative preconditioned approach reduces memory requirements compared to direct methods.
		\end{abstract}
	\end{frontmatter}
	\section{Introduction}
	
	Solving Maxwell's equations is one of the major problems of mathematical modeling. Of particular interest is calculating harmonic electromagnetic fields in geometrically complex domains. In this case, as a rule, the vector finite element method with basis functions from the Nedelec functional space is used~\cite{Nedelec1980,Nedelec1986}. The thus obtained discrete problem gives a physically correct approximation of Maxwell's equations providing, in a natural way, the discontinuity of normal components and the continuity of tangential components of the electromagnetic field without adding fictitious resonance frequencies to the spectrum. We assume a time dependence of the form $e^{\ii \omega t}$.
	
	Let $\Omega \subset \mathbb{R}^d,\ d=2 \text{ or }3$  be a bounded Lipschitz domain with two disjoint connected boundaries~$\Gamma$ and $S$ with outward unit normal  $\mathbf{n}$. By solving the time-harmonic Maxwell equations under perfect conducting boundary conditions on $S$ and impedance conditions~on~$\Gamma$:
	\begin{subequations}
		\begin{empheq}{align}
			\nabla \times \left(\mu_r^{-1} \nabla \times \tilde{\textbf{E}}\right) - \kappa^2\epsilon_r^c\tilde{\textbf{E}}=\tilde{\textbf{F}}                                          \; \text{in}\; \Omega, \\
			\textbf{n} \times \tilde{\textbf{E}}=0                                                                                                                                 \; \text{on}\;S,       \\
			\mu_r^{-1}\left(\nabla \times \tilde{\textbf{E}}\right)\times \textbf{n}-\ii \kappa \lambda\left(\textbf{n}\times  \tilde{\textbf{E}}\right)\times \textbf{n}=0  \; \text{on}\; \Gamma,
		\end{empheq}
		\label{HelmVector}
	\end{subequations}
	we can get a time-harmonic electric field  $\tilde{\textbf{E}}:\Omega \rightarrow \mathbb{C}^d$ corresponding to a given current density $\tilde{\mathbf{F}} \in L_2(\mathbb{R})^{d}$.
	The complex absolute and relative permittivity are defined as:
	$$
	\epsilon_c(\mathbf{x}) = \epsilon(\mathbf{x}) + \frac{\ii\sigma(\mathbf{x})}{\omega}, \quad \epsilon_r^c(\mathbf{x}) = \frac{\epsilon_c(\mathbf{x})}{\epsilon_0}.
	$$
	The wave number is related to the frequency $\omega$ and material parameters by:
	$$
	\kappa = \omega\sqrt{\mu_0\epsilon_0\mu_r(\mathbf{x})\epsilon_r^c(\mathbf{x})}, \quad \mu_r=\frac{\mu(\mathbf{x})}{\mu_0},
	$$ 
	where $\epsilon,\mu$ are the electric and magnetic permeability (for a vacuum $\epsilon_0,\mu_0$), $\sigma$ is the conductivity, $\lambda$ is the impedance coefficient. \fix{The proposed method is applicable for $\sigma \geq 0$, however, the impedance boundary condition is considered for the lossless case $\sigma=0$.}
	
	The Hilbert spaces corresponding to the Maxwell equation are defined by 
	$$
	H\left(\mathrm{curl};\Omega\right)=\left\{\textbf{u}\in \left(L_2\left(\Omega\right)\right)^d:\nabla\times \textbf{u} \in \left(L_2\left(\Omega\right)\right)^d \right\},
	$$
	$$
	H_0=\left\{ \mathbf{u}\in H(\mathrm{curl};\Omega)\;:\; \mathbf{n}\times\mathbf{u}=0\ \text{on }S\right\},
	$$
	equipped with the following curl norms:
	$$
	\|\textbf{u}\|_{H\left(\mathrm{curl};\Omega\right)}^2={\|\textbf{u}\|^2_{\left(L_2\left(\Omega\right)\right)^d}+\|\nabla\times \textbf{u}\|^2_{\left(L_2\left(\Omega\right)\right)^d}},
	$$
	$$
	\fix{\|\mathbf{u}\|^2_{H_0}=\|\mathbf{u}\|^2_{H\left(\mathrm{curl};\Omega\right)}+\|\mathbf{n}\times \mathbf{u}\|^2_{\left(L_2\left(\Gamma\right)\right)^{d-1}}}.
	$$
	To simplify the notation, we introduced the following abbreviations for inner products
	and boundary integrals:
	$$
	\left(\mathbf{u},\mathbf{{v}}\right)_\Omega=\int_{\Omega}\mathbf{u}\cdot{\mathbf{\bar{v}}}d\Omega, \quad \left<\mathbf{u},\mathbf{{v}}\right>_{\Gamma}=\int_{\Gamma}\mathbf{u}\cdot{\mathbf{\bar{v}}}d\Gamma,
	$$
	where the overbar denotes complex conjugation.
	Using the Galerkin method, we arrive at the variational problem of finding $\tilde{\textbf{E}} \in H_0$ such that
	\begin{equation}
		\label{variational_maxwell}
		\left(\mu_r^{-1}\nabla \times \tilde{\textbf{E}},\nabla \times \mathbf{v} \right)_\Omega-\kappa^2\epsilon_r^c\left(\tilde{\textbf{E}},\mathbf{v} \right)_\Omega-\ii\kappa\lambda\left<\tilde{\textbf{E}}_T,\mathbf{v}_T\right>_{\Gamma}=
		\left(\tilde{\mathbf{F}},{\mathbf{v}}\right)_\Omega, \quad \forall\; {\mathbf{v}}\in {H_0}.
	\end{equation}
	If the calculation domain $\Omega$ is approximated, for example, by a triangular $(d=2)$ or tetrahedral $(d=3)$ grid, and the space $H_0$ is approximated by a finite-dimensional Nedelec subspace $H_0^h \in H_0$~\cite{Monk1993}, the solution of the variation problem (\ref{variational_maxwell})  can be obtained by solving the system of linear algebraic equations
	\begin{equation}
		A \mathbf{X}=\mathbf{Y}, \quad A \in \mathbb{C}^{n\times n},\quad \mathbf{X},\mathbf{Y}\in \mathbb{C}^n.
		\label{main_linear_system}
	\end{equation}
	
	Significant computational difficulties arise when solving an ill-conditioned system of linear algebraic equations (SLAE) of the form (\ref{main_linear_system}). Although many computational methods of linear algebra have been developed for this purpose, the problem of solving high-order indefinite SLAEs has not yet been fully solved, especially with iterative methods~\cite{Ernst2012}. In engineering calculations, preference is given to direct algorithms of solving SLAEs with various variants of the Gaussian elimination~\cite{duff2017}, which is due to the fact that they are computationally reliable, but fast convergence of iterative methods is not guaranteed. However, the direct methods require great numbers of arithmetic operations, large RAM, and use of multiprocessor computing systems even for solving relatively small problems.
	
	An alternative to direct methods is iterative algorithms \cite{VanDerVorst_2003} which, for many problems of mathematical modeling, make it possible to obtain SLAE solutions with fewer arithmetic operations. For the vector Helmholtz equation with a sign-definite operator, various approaches have been developed, which are mainly based on algebraic multigrid methods \cite{Hiptmair1998,Jones2006,Hiptmair2007,Vassilevski2009}. In some particular cases, multi-grid methods can be used to solve sign-indefinite problems \cite{NECHAEV2008}. A geometric multigrid method can be used when it is possible to construct a sequence of nested calculation grids \cite{Arnold2000}, which is difficult in geometrically complex domains.
	
	In the general case, the difficulty in solving the ill-conditioned system of equations (\ref{main_linear_system}) is that the principal part of (\ref{HelmVector}a) has a nontrivial null-space, since \mbox{$\nabla \times \left(\mu_r^{-1}\nabla \times \left(\nabla \psi \right)\right) = 0$} for any three times differentiable scalar function $\psi$. Here the dimensionality of the null-space is proportional to the number of grid nodes \cite{NullSpace_Maxwell}, which does not allow using iterative algorithms for solving the SLAEs efficiently. The sign-indefiniteness and the presence of a nontrivial null-space of the matrix of a SLAE do not prevent the use of iterative methods; however, their combination when solving Maxwell's equations greatly reduces the convergence rate. \fix{It should be noted that the null-space is associated not only with the gradient of scalar functions but also with the improper approximation of the divergence operator in the finite element method.}
	
	To relax the requirements for computing resources, alternative approaches based on the principle of maximum amplitude have been developed \cite{Tikhonov1948,Tikhonov1990}. For example, for the scalar Helmholtz equation	
	\begin{equation}
		\Delta \tilde{u}(\textbf{x})+k^2\tilde{u}(\textbf{x})=-\tilde{f}(\textbf{x}), \quad k^2=\frac{\omega^2}{v^2}, \quad  \mathbf{x} \in \mathbb{R}^d, \; d=2\text{ or }3,
		\label{helm}
	\end{equation}
	where $k^2$ is the squared wave vector modulus, $\omega$ is the cyclic frequency, and  $v$ is the wave propagation speed in the medium. The solution to equation~(\ref{helm}) is understood as the limit
	\begin{equation}
		\tilde{u}(\textbf{x})=\lim\limits_{t\to \infty}u(\textbf{x},t)e^{\mp \ii \omega t}, \quad \ii=\sqrt{-1}. 
		\label{limit}
	\end{equation}
	\fix{(Requires radiation conditions; invalid for bounded domains without absorbing boundary conditions.)}
	Here the function $u(\textbf{x},t)$ is a solution to the following Cauchy problem:
	\begin{equation}
		\begin{array}{l}
			\Delta u(\textbf{x})-\frac{1}{v^2}\partial^2_{tt}u(\textbf{x})=-\tilde{f}(\textbf{x})e^{\pm \ii \omega t},
			\\\\ u(\textbf{x},0)=\partial_{t} u(\textbf{x},0)=0.
		\end{array}
		\label{Wave_eq}
	\end{equation}
	Along with the limiting absorption principle or the Sommerfeld radiation conditions~\cite{Sommerfeld1912}, the limiting amplitude principle makes it possible to identify the unique solution that is physically correct. A necessary condition of correctness of the representation (\ref{limit})  is that the function $\tilde{f}(\mathbf{x})$ must be orthogonal to all eigenfunctions of the steady problem (\ref{helm}). Hence, for bounded areas the limiting amplitude principle is not valid. From the computational point of view, no need to solve ill-conditioned sign-indefinite SLAEs is an advantage of this class of methods. A shortcoming is that, due to the low convergence rate, which depends on the model parameters and the dimensions of the calculation domain, it may be necessary to calculate a significant number of wave field periods in order to achieve a given accuracy by means of "pseudo transient” computations. 
	
	In order to reduce the computational costs, some advanced approaches known as "controllability methods"{} have been proposed. They require multiple solutions of Maxwell’s equations in time to minimize a penalty functional~\cite{Bristeau1998,Grote2019,ChaumontFrelet2022}, whose minimum is reached on a solution that is close to that of equation (\ref{HelmVector}). In contrast to the limiting amplitude principle, the use of the additional functional reduces the number of problems to be solved in time; however, the thus obtained solution should be understood in a very weak sense. The main difficulty of this class of methods is that a solution to the functional may be nonunique. In addition, there is no guarantee that the solution providing a minimum to the functional will be close to the solution of the problem (\ref{HelmVector}). Papers~\cite{Appel2020,Peng2022}  consider various wave field filtering procedures to identify the unique solution obtained in this way.
	
	When calculating a large number of wavefield periods, small grid steps in space and time must be used to ensure sufficient accuracy. For explicit schemes, the integration time step is limited by the Courant stability condition~\cite{Courant1928}. Implicit schemes allow calculations with larger time steps. This may cause parasitic non-physical oscillations in the wave field, and makes it necessary to solve sign-indefinite SLAEs.
	
	The present study deals with preconditioning of a SLAE of the form (\ref{main_linear_system}). We construct a preconditioning operator based on repeatedly solving Maxwell’s equations in time using Laguerre transform. This makes it possible to obtain a series of problems with a sign-definite operator, for which the discrete problem can be effectively solved by multigrid algorithms. The use of Laguerre transform eliminates numerical dispersion in time, and there is no need to perform calculations with very fine time resolution. In contrast to "controllability methods"{} and other similar approaches, solving the problem (\ref{main_linear_system}) directly does not require any additional auxiliary functionals, which makes it easier to analyze the accuracy of the results.
	
	\section{Laguerre based preconditioner}
	Consider integral Laguerre transform in time for preconditioning of SLAE~(\ref{main_linear_system}). This transform has been used to solve direct and inverse problems in modeling of seismic and electromagnetic wave fields \cite{Mikhailenko2008,Mastryukov2013,fatab2011,Terekhov2015206,Mikhailenko1999}. It has been shown that with this approach no sign-indefinite SLAEs need be solved and spectral accuracy in approximating the time derivative is provided \fix{\cite{Mikhailenko1999,Terekhov2018}}. Also, Laguerre transform makes it possible to stabilize unstable difference schemes when solving problems of wave field continuation from the surface into depth \cite{Terekhov2017,Terekhov2018}
	
	\subsection{Laguerre transform}
	Consider Laguerre functions \cite{NIST:DLMF}, which are defined as
	\begin{equation*}
		l_m(t)=e^{-t/2}\mathsf{L}_m(t),\quad t\geq 0,\; m=0,1,2...,
		\label{laguerre_function}
	\end{equation*}
	where $\mathsf{L}_m(t)$ is the Laguerre polynomial of degree $m$, which is defined by the Rodrigues formula
	$$
	\mathsf{L}_m(t)=\frac{ e^{t}}{m!}\frac{d^m}{dt^m}\left(t^me^{-t}\right)=\frac{1}{m!}\left(\frac{d}{dt}-1\right)^mt^m.
	$$
	We will use $L_2[0,\infty)$ to denote the space of square integrable functions \mbox{$f:[0,\infty)\rightarrow \mathbb{R}$}
	\begin{equation*}
		L_2[0,\infty)=\left\{f:\int_{0}^{\infty}|f(t)|^2dt< \infty\right\}.
	\end{equation*}
	The Laguerre functions are a complete orthonormal system in $ L_2[0,\infty)$
	\begin{equation*}
		\int_{0}^{\infty}l_m(t)l_n(t)dt=
		\left\{\begin{array}{ll}
			0, & m \neq n , \\
			1, & m=n.
		\end{array}\right.
		\label{lag_orho}
	\end{equation*}
	This guarantees that for any function $f(t)\in L_2[0,\infty)$ there is a Laguerre expansion
	\begin{equation}
		\label{series_lag}
		f(t)\sim\eta \sum_{m=0}^{\infty}\bar{a}_m l_m(\eta t), \quad \bar{a}_m=\int_{0}^{\infty}f(t)l_m(\eta t)dt, \quad t  \geq 0,
	\end{equation}
	where $\eta>0$ is a scaling parameter for the Laguerre functions to increase the convergence rate of the series \cite{Weeks1966}.
	
	Let us consider bilateral Laguerre transform for functions  $g(t) \in L_2(\mathbb{R})$~\cite{Keilson1981}:
	\begin{equation}
		\begin{array}{ll}
			\displaystyle	g(t)\sim\eta \sum_{m=-\infty}^{+\infty}\bar{b}_m \phi_m(\eta t),\quad \bar{b}_m=\int_{-\infty}^{+\infty}g(t)\phi_m(\eta t)dt,
		\end{array}
		\label{series_lag_R}
	\end{equation}
	where  $\phi_m(t)$  are functions that are orthonormal on the entire real axis and defined as follows:
	\begin{equation}
		\phi_m(t)=\left\{
		\begin{array}{lc}
			\left.\begin{array}{ll}
				l_m(t), & t \ge 0 \\
				0,      & t<0
			\end{array}\right\}\quad m=0,1,2,...\quad , \\
			\left.\begin{array}{ll}
				-l_{-m-1}(-t), & t < 0  \\
				0,             & t\ge 0
			\end{array}\right\}\quad m=-1,-2,...\quad.
		\end{array}\right.
	\end{equation}
	Then Plancherel's extension of the Fourier transform for the function $\phi_m(t)$ has the form \cite{Weber1980}:
	\begin{equation}
		\tilde{\phi}_m(\omega)=\mathcal{F}\{\phi_m(\eta t)\}=\frac{1}{2\pi}\int_{-\infty}^{+\infty} \phi_m(\eta t)e^{-\ii\omega t}dt=\frac{1}{2\pi}
		\frac{\left(\ii \omega-\eta/2\right)^m}{\left(\ii \omega+\eta/2\right)^{m+1}}, \quad m \in \mathbb{Z}.
	\end{equation} 
	Hence, the Fourier transform is expressed in terms of the Laguerre series coefficients as follows:
	\begin{equation}
		\label{Laguerre_Fourier}
		\tilde{g}(\omega)=\mathcal{F}\{g(t)\}\sim\eta\sum_{m=-\infty}^{+\infty}\bar{b}_m\tilde{\phi}_m(\omega).
	\end{equation}
	
	This relation will be used to construct a preconditioning procedure for solving SLAEs of the form (\ref{main_linear_system}). The calculation time can be decreased by adjusting the number of terms in the Laguerre series expansion, since for preconditioning the solution is allowed to be calculated approximately.
	
	\subsection{The Laguerre transform for Maxwell's equation}
	Taking into account the operator properties of Fourier transform \mbox{$i\omega \leftrightarrow \partial_t$} and \mbox{$-\omega^2 \leftrightarrow \partial_{tt}$}, we can write the problem (\ref{HelmVector}) in time as follows:
	\begin{equation}
		\begin{array}{rl}
			\nabla \times \left(\mu_r^{-1} \nabla \times \textbf{E}\right) 
			\fix{+ \mu_0 \mu_r \left( \epsilon_r \epsilon_0 \partial_{tt}\mathbf{E} + \sigma \partial_t\mathbf{E} \right) = \mathbf{F}} & \text{in}\; \Omega,\\
			\textbf{n} \times \textbf{E}=0                                                                                                                                 & \text{on}\;S,\\
			\mu_r^{-1}\left(\nabla \times {\textbf{E}}\right)\times \textbf{n}-\lambda\sqrt{\mu_r \epsilon_r}\partial_t\left(\textbf{n}\times  {\textbf{E}}\right)\times \textbf{n}=0 & \text{on}\; \Gamma \fix{\;(\sigma=0)}.
			\label{time_maxwell}
		\end{array}
	\end{equation}
	A solution to equation (\ref{time_maxwell}) will be sought for in the form of a series of Laguerre functions:
	\begin{equation}
		\mathbf{E}(\mathbf{x},t)\approx\eta \sum_{m=0}^{M-1}\mathbf{\bar{E}}_ml_m(\eta t),
		\label{solution_approx}
	\end{equation}
	where the number of expansion coefficients is chosen according to the required accuracy of approximation of the solution on an interval, $t\in [0,t_1]$. Taking into account $f(0)=f'(t)=0$ and assuming \mbox{$\lim_{t\rightarrow \infty}f(t)=\lim_{t\rightarrow \infty}f'(t)=0$}, one can show \cite{Integral_Transform} that
	\begin{equation}
		\int_{0}^{\infty}\frac{d^{k}f(t)}{d t^k}l_m(\eta t)dt=\left(\frac{\eta}{2}\right)^{k}\bar{f}_m+\Phi_k(\bar{f}_m), \quad m=0,1,2...,
		\label{partial_t_lag}
	\end{equation}
	where for $k=1,2$ we have
	\begin{equation*}
		\fix{\Phi_1(\bar{f}_m)\equiv \eta\sum_{k=0}^{m-1}\bar{f}_k, \quad\Phi_2(\bar{f}_m)\equiv\eta^2\sum_{k=0}^{m-1}(m-k)\bar{f}_k, \quad m=0,1,2....}
		\label{series_prop}
	\end{equation*}
	Then, assuming  $\textbf{E}|_{t=0}=\partial_t \textbf{E}|_{t=0}= \lim_{t\rightarrow \infty}\mathbf{E}(t)=\lim_{t\rightarrow \infty}\partial_t\mathbf{E}(t)=0$ and applying the Laguerre transform to the system~(\ref{time_maxwell}), taking into account the properties~(\ref{partial_t_lag}), we obtain
	\begin{equation}
		\begin{array}{rl}
			\nabla \times \left(\mu_r^{-1} \nabla \times \bar{\mathbf{E}}_m\right)+
			\fix{\beta_1 \bar{\mathbf{E}}_m=\mathbf{\bar{F}}_m+\sigma\mu_r/\epsilon_0\Phi_1(\bar{\mathbf{E}}_m)+\mu_r\epsilon_r\Phi_2(\bar{\mathbf{E}}_m)} & \text{in}\; \Omega,\\
			\mathbf{n} \times \bar{\textbf{E}}_m=0                                                                                                                                 & \text{on}\;S,\\
			\mu_r^{-1}\left(\nabla \times {\mathbf{\bar E}}_m\right)\times \mathbf{n}-\beta_2\left(\mathbf{n}\times{\mathbf{\bar{E}}_m}\right)\times \mathbf{n}=\lambda\sqrt{\mu_r \epsilon_r}\Phi_1\left(\left(\mathbf{n}\times  {\mathbf{\bar{E}}_m}\right)\times \mathbf{n}\right) & \text{on} \; \Gamma,
		\end{array}
		\label{Maxwell_laguerre}
	\end{equation}
	where \fix{$\beta_1(\mathbf{x}) = \mu_0 \mu_r(\mathbf{x}) \left( \epsilon_r(\mathbf{x}) \epsilon_0 \frac{\eta^2}{4} + \sigma(\mathbf{x}) \frac{\eta}{2} \right)$, $\beta_2(\mathbf{x})=\frac{\eta}{2}\lambda\sqrt{\mu_r(\mathbf{x}) \epsilon_r(\mathbf{x})}$}.
	
	To calculate the expansion coefficients of the series (\ref{solution_approx}), we approximate the problem (\ref{Maxwell_laguerre}) by using a vector finite element method for the same subspace $H^h_0$  and the grid as in the problem (\ref{HelmVector}):
	\begin{equation}
		\label{variational_maxwell_precond}
		\left(\mu_r^{-1}\nabla \times \bar{\textbf{E}}_m,\nabla \times \mathbf{v} \right)_\Omega \fix{+} \beta_1\left(\bar{\textbf{E}}_m,\mathbf{v} \right)_\Omega-\beta_2\left<\bar{\textbf{E}}_{T,m},\mathbf{v}_T\right>_{\Gamma}=
		\left(\bar{\mathbf{Q}}_m,{\mathbf{v}}\right)_\Omega, \quad \forall\; {\mathbf{v}}\in {H_0}.
	\end{equation}
	In what follows, we assume that all material parameters $\epsilon(\mathbf{x}),\mu(\mathbf{x}),\sigma(\mathbf{x})$ and the impedance coefficient $\lambda$ are real-valued, and the initial time-domain data are real. Under these assumptions, the preconditioning matrix $B$ obtained after Laguerre expansion and finite element approximation on a real basis is real-valued.
	Then the following SLAE is solved several times with the same matrix and different right-hand sides:
	\begin{equation}
		B\mathbf{\bar E}_m=\mathbf{\bar G}_{m-1}+\mathbf{\bar{F}}_m, \quad B \in \mathbb{R}^{n\times n},\quad \mathbf{\bar E}_m,\mathbf{\bar G}_{m-1},\mathbf{\bar{F}}_m\in \mathbb{R}^n
		\label{preconditioner}
	\end{equation}
	where the vector  $\mathbf{\bar G}_{m-1}$ is a linear combination of   $\mathbf{\bar E}_k$, $k=0,1,2,...,m-1$ and \mbox{$\mathbf{\bar E}_{-1}\equiv0$}.
	At \mbox{$\beta_1>0$} (i.e. $\eta > \frac{2\sigma(\mathbf{x})}{\epsilon(\mathbf{x})}, \forall \mathbf{x} \in \Omega$), this sign-definite SLAE can be effectively solved by an algebraic multigrid algorithm \cite{Vassilevski2009}. Compared to direct methods for solving SLAEs, multigrid methods need much less RAM.
	
	\section{Choosing the right-hand side}
	Consider how to choose the right-hand side for the problem (\ref{variational_maxwell_precond}),(\ref{preconditioner})  to precondition SLAE~(\ref{main_linear_system}).
	Since the functions $\sin(\omega t),\cos(\omega t)$ do not belong to the space $L_2[0,\infty]$, the Laguerre series will diverge for the right-hand sides of the form $\mathbf{F}(\mathbf{x},t)=\mathbf{\tilde{F}}(\textbf{x})e^{\ii\omega t}$. However, convergence can be achieved if the right-hand side is multiplied by a function defining a rectangular window:
	$$
	\Pi(t)=\begin{cases}1,& 0\le t\le\tau,\\ 0,&\text{otherwise.}\end{cases}
	$$
	Hence,
	$$\mathbf{F}(\mathbf{x},t)=\mathbf{\tilde{F}}(\textbf{x})e^{\ii\omega t}\Pi(t).$$
	Then the harmonic electromagnetic field will be approximately calculated as follows:
	\begin{equation}
		\fix{\mathbf{\tilde{E}}(\mathbf{x},\omega) \approx \eta \sum_{m=0}^{M} \mathbf{\bar{E}}_m(\mathbf{x}) \cdot \tilde{\phi}_m(\omega).}
		\label{Approximate_E}
	\end{equation}
	Since the function $\Pi(t)$ is not differentiable, the Laguerre series (\ref{series_lag}) will be slowly converging. On the other hand, as the duration of the signal decreases with time, the width of its spectrum increases, since after a harmonic source is multiplied by the function $\Pi(t)$ the spectrum is the convolution of the delta function with a $sinc$ function, which causes the Gibbs phenomenon. Also, the "spectral leakage"{} effect makes it a smooth function of $\omega$, which increases the convergence rate of the series (\ref{Approximate_E}). Thus, one cannot provide equally fast convergence of the Laguerre series to approximate the source in the time domain and to approximate the electromagnetic field in the spectral domain.
	
	Window functions such as Kaiser, Hanning, and Chebyshev ones have been additionally tested. However, at this stage of the study their use does not provide the convergence of the iterative process for solving the SLAE. This can be explained by the fact that the best means for retaining such characteristics as variance compensation, dispersion factors, total energy, and coherent gain is a rectangular window~\cite{Prabhu2018}, which turned out to be important in preconditioning.
	\section{Numerical experiments}
\begin{figure}[!htb]
		\centering
		\subfloat{
			\includegraphics[width=0.9\textwidth]{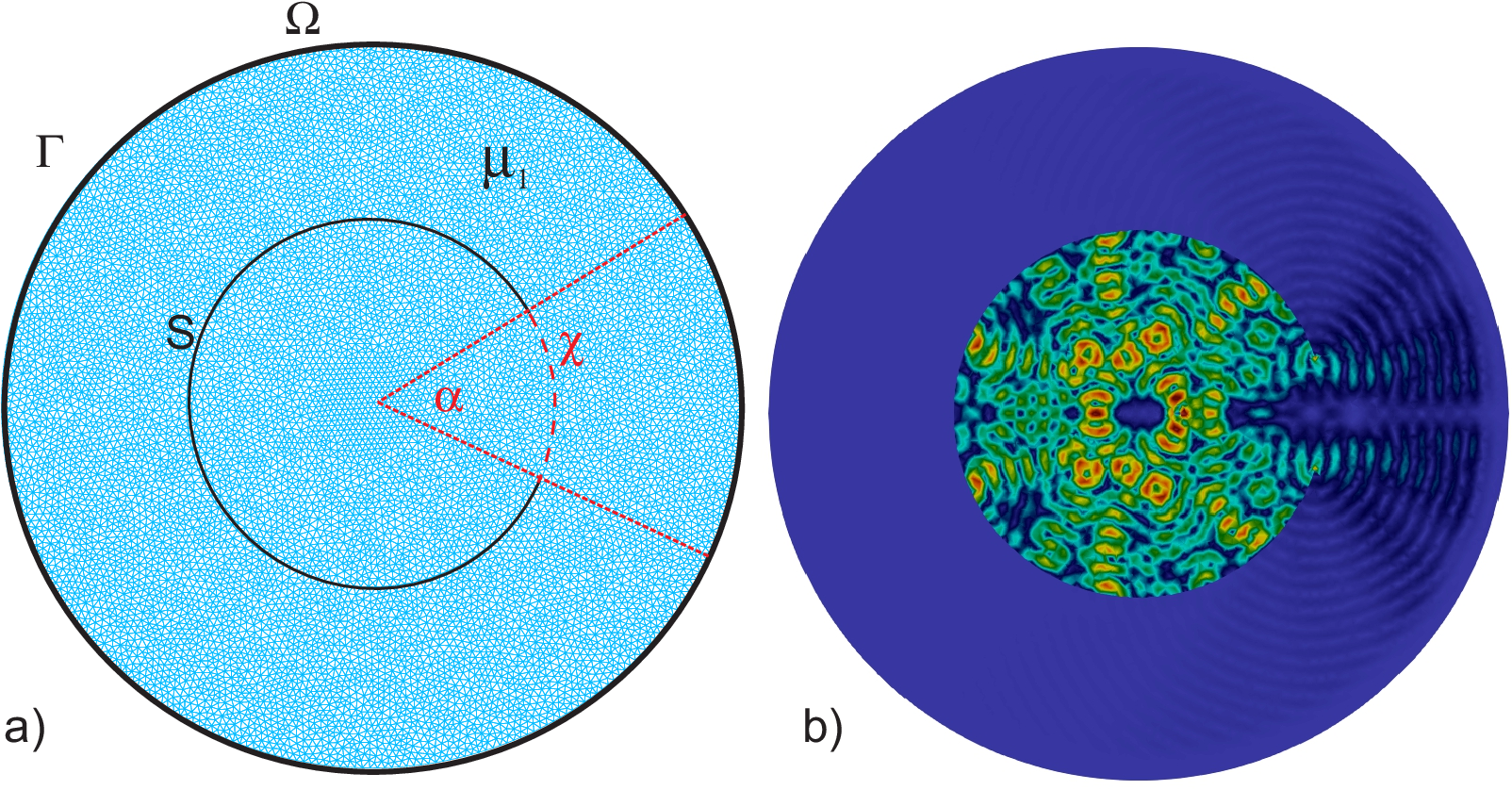}}
		\caption{ a) Homogeneous model of medium with internal PEC boundary $S$, b) wave field magnitude for a source located at the center of the calculation domain}
		\label{pic:test_model1}	
	\end{figure}

\begin{figure}[!htb]
		\centering
		\includegraphics[width=0.9\textwidth]{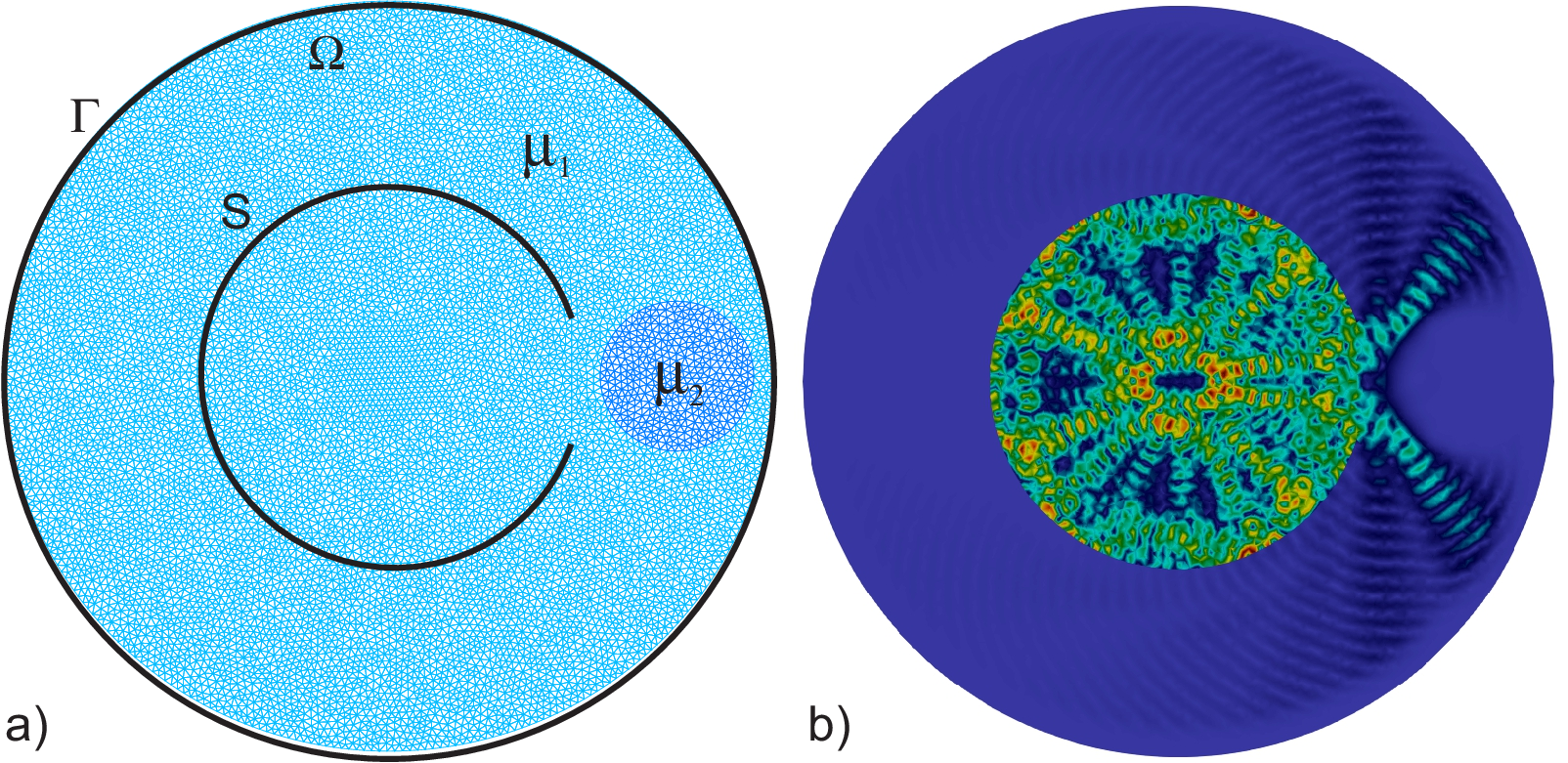}
		\caption{a)~Heterogeneous model of the medium with internal PEC boundary $S$, b)~wavefield magnitude for a source with frequency $2\pi/100$ located at the center of the calculation domain}
		\label{pic:test_model2}
	\end{figure}

	Consider solving equation  (\ref{HelmVector})  for the two media models shown in Figures \ref{pic:test_model1}a and \ref{pic:test_model2}a. The parameters of the medium are as follows: permeability $\mu_1=30$, $\mu_2=5$, permittivity $\epsilon_r=1$. A point source with frequency  $\omega=2\pi/100$ is located at the center of the calculation domain $\Omega$. After a finite-element approximation of the variational problem (\ref{variational_maxwell}) with second-order Nedelec's elements, we solve a SLAE of the form (\ref{main_linear_system}), for which we use the GMRES($k$) method~\cite{Saad,Saad1986} and the preconditioning procedure (\ref{solution_approx}), (\ref{preconditioner}) with the Laguerre transform parameter $\eta=0.15$ and time window duration $\tau=4\pi/\omega$. The coefficients $\bar{\mathbf E}_m$ of the series (\ref{Laguerre_Fourier}) will be calculated with accuracy $\varepsilon_{lag}=10^{-5}$. \fix{The parameter $\eta$ was chosen to balance spectral localization and time support.} The scaling parameter $\eta$ and window length $\tau$ were chosen to balance spectral concentration and time support; $\varepsilon_{lag}$ is the stopping tolerance for solving the Laguerre-projected systems. A short sensitivity study (varying $\eta$ and $M$) is recommended to justify these choices. As has been shown by computational experiments, to solve the problem (\ref{preconditioner}) we need only a few iterations of the conjugate gradient method with a multigrid preconditioner called the Auxiliary-space Maxwell Solver \cite{Vassilevski2009}.
	
\begin{figure}[!htb]
		\centering
		\includegraphics[width=1\textwidth]{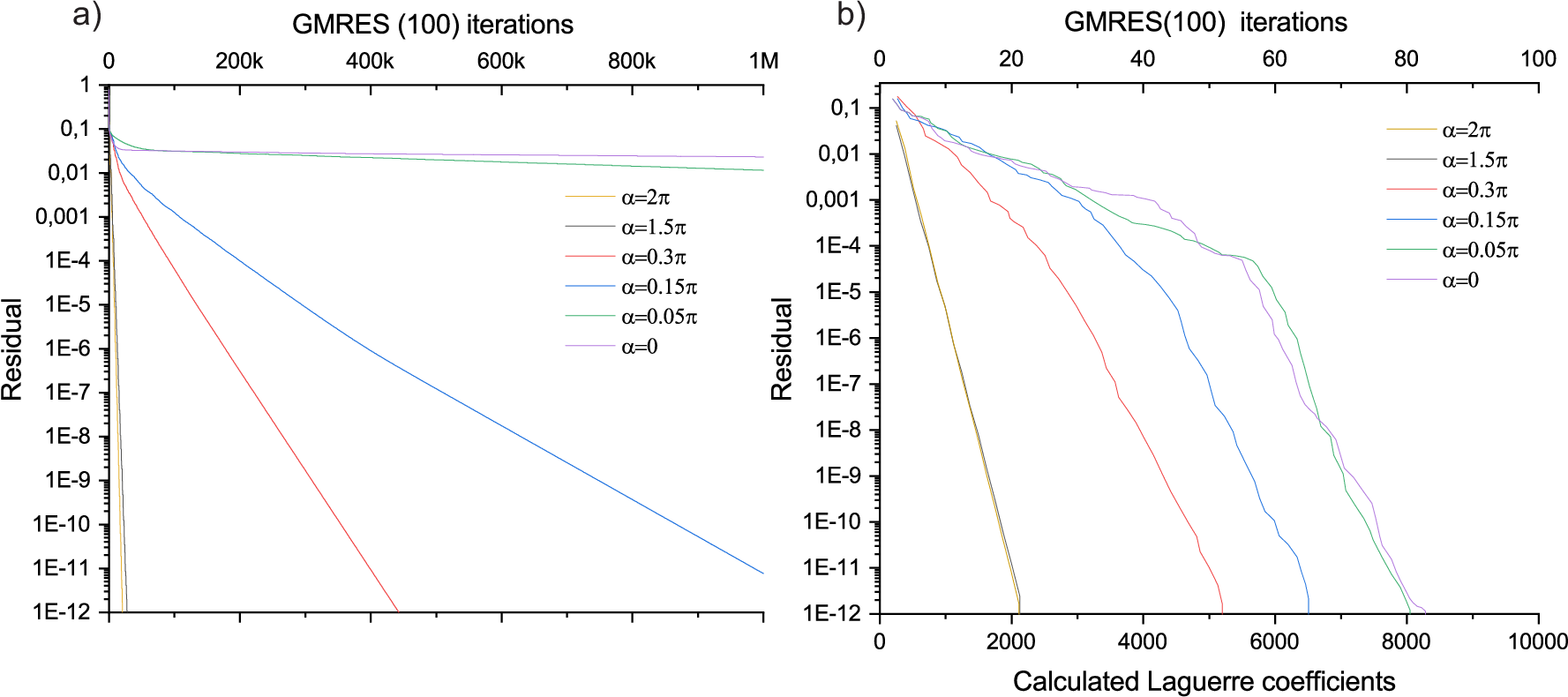}
		\caption{Residual for the GMRES($k$) method versus the number of iterations for the homogeneous model (Fig.~\ref{pic:test_model1}a): a) without preconditioning and b) with preconditioning}
		\label{pic:residual1}
	\end{figure}

	Let us investigate the convergence of the GMRES($k$) method, $k=100$, without preconditioning. The wave field magnitude for the first and second models is shown in Figs.~\ref{pic:test_model1}b and~\ref{pic:test_model2}b, respectively. As can be seen in Fig.~\ref{pic:residual1}a, the convergence rate depends on the length $\chi(\alpha)$ of the segment PEC of the boundary $S$. If the boundary $S$ is absent ($\alpha=2\pi$) and there are no reflected waves, the iterative process converges very fast. If the boundary is almost closed or completely closed ($\alpha\leq0.05\pi$), the GMRES($k$) method stagnates.
	
\begin{figure}[!htb]
		\centering
		\includegraphics[width=1\textwidth]{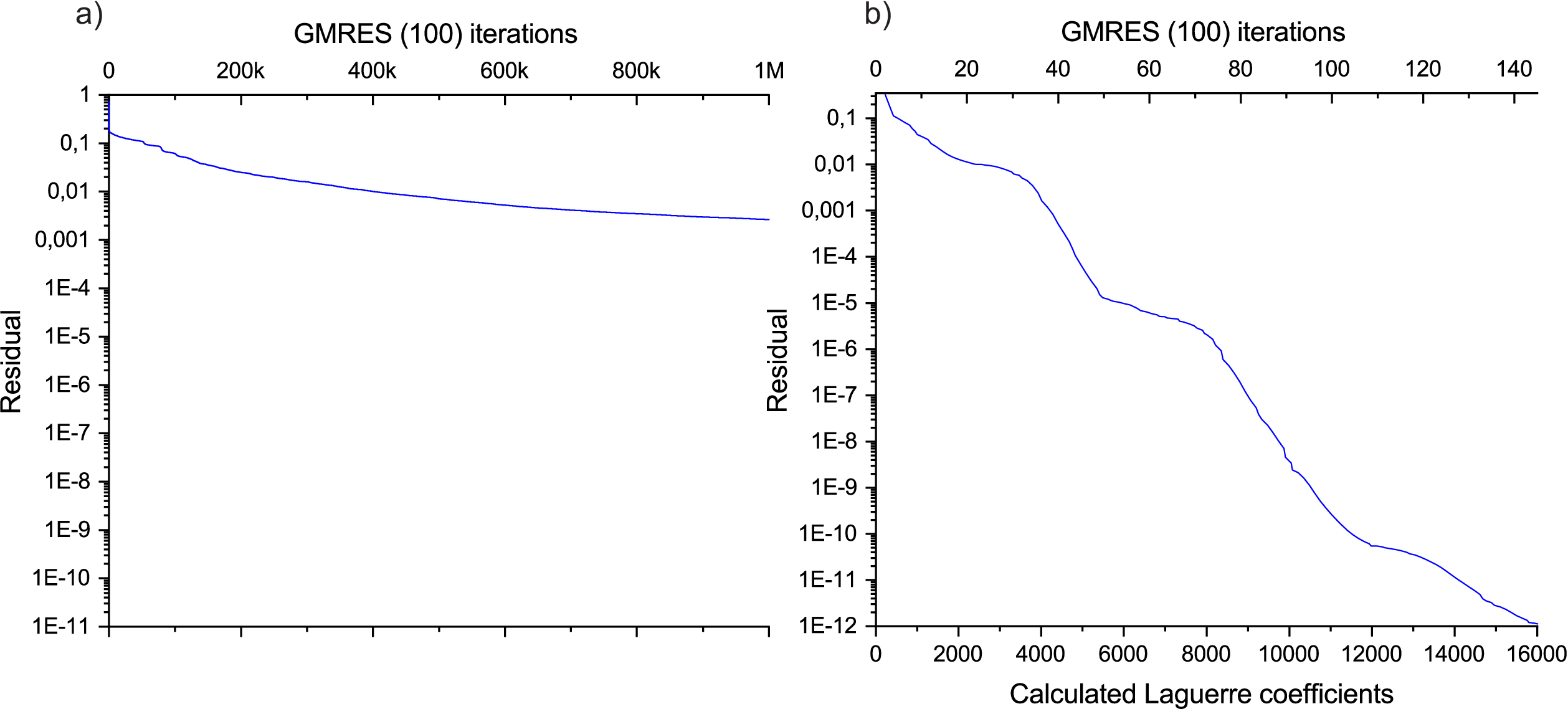}
		\caption{Residual for the GMRES($k$) method versus the number of iterations for the heterogeneous model (Fig.~\ref{pic:test_model2}a): a) without preconditioning and b) with preconditioning}
		\label{pic:residual2}
	\end{figure}

	For the heterogeneous model stagnation is observed (Fig.~\ref{pic:residual2}a) already at $\alpha=0.05\pi$, although for the homogeneous model with the same segment length $\chi(\alpha)$ the iterative process converges. Thus, without preconditioning the iterative methods converge very slowly for heterogeneous models, in particular for PEC boundary conditions resulting in multiple internal reflections. By contrast, the GMRES ($k$) method with a preconditioner of the form (\ref{solution_approx}), (\ref{preconditioner})  provides convergence for both the homogeneous model and the heterogeneous one. As the slot length $\chi(\alpha)$ decreases, there appear multiple reflections of waves from the boundary $S$; and the waves need more time to leave the calculation domain $\Omega$. As a consequence, the convergence rate for small values of  $\alpha$ decreases, since the preconditioning operator approximates the solution in some time interval $[0,\tau]$, which does not guarantee that the solution of the problem (\ref{preconditioner}) is close to that of the problem (\ref{main_linear_system}).

	\section{Conclusions}
	This paper considered solving Maxwell's equations in the frequency domain. The sign-indefiniteness, poor conditioning of the SLAE, and the presence of a nontrivial null-space make it difficult to reach high convergence rates of iterative methods when considering many engineering problems. The algorithms based on the Gaussian elimination require a significant number of arithmetic operations, as well as large RAM, even in solving relatively small problems. The main requirement in the development of the present approach was to decrease the minimum RAM needed. For this, a preconditioning procedure was proposed, which is based on integral Laguerre transform in time. This makes it possible to obtain a sign-definite SLAE, which can be solved by multigrid methods that do not require large RAM.
	
	In the two-dimensional case the approach is not efficient, since it needs more calculation time than the direct methods. However, in the three-dimensional case the efficiency will be much better, since the RAM required is several orders of magnitude less than that of the direct methods. The calculation time will also decrease, since the computational complexity of multigrid methods is less than that of direct algorithms. Also note that the computational procedure under consideration does not stagnate, which is due to the filtration properties of the Laguerre transform which were used to suppress the unstable solution components in our previous studies.
	\newpage
	\bibliography{base}

\begin{thebibliography}{10}

\bibitem{Nedelec1980}
J.C. Nédélec.
\newblock Mixed finite elements in {$\mathbb{R}^3$}.
\newblock {\em Numerische Mathematik}, 35(3):315--341, 1980.

\bibitem{Nedelec1986}
J.C. Nédélec.
\newblock A new family of mixed finite elements in {$\mathbb{R}^3$}.
\newblock {\em Numerische Mathematik}, 50:57--81, 1986.

\bibitem{Monk1993}
P.~Monk.
\newblock An analysis of n{\'{e}}d{\'{e}}lec{\textquotesingle}s method for the
  spatial discretization of maxwell{\textquotesingle}s equations.
\newblock {\em Journal of Computational and Applied Mathematics},
  47(1):101--121, June 1993.

\bibitem{Ernst2012}
O.~G. Ernst and M.~J. Gander.
\newblock {\em Why it is Difficult to Solve Helmholtz Problems with Classical
  Iterative Methods}, pages 325--363.
\newblock Springer Berlin Heidelberg, Berlin, Heidelberg, 2012.

\bibitem{duff2017}
I.S. Duff, A.M. Erisman, and J.K. Reid.
\newblock {\em Direct Methods for Sparse Matrices}.
\newblock Oxford University Press, January 2017.

\bibitem{VanDerVorst_2003}
H.~A. van~der Vorst.
\newblock {\em Iterative Krylov Methods for Large Linear Systems}, volume~13.
\newblock Cambridge University Press, 2003.

\bibitem{Hiptmair1998}
R.~Hiptmair.
\newblock Multigrid method for maxwell{\textquotesingle}s equations.
\newblock {\em {SIAM} Journal on Numerical Analysis}, 36(1):204--225, January
  1998.

\bibitem{Jones2006}
J.~Jones and B.~Lee.
\newblock A multigrid method for variable coefficient
  maxwell{\textquotesingle}s equations.
\newblock {\em {SIAM} Journal on Scientific Computing}, 27(5):1689--1708,
  January 2006.

\bibitem{Hiptmair2007}
R.~Hiptmair and J.~Xu.
\newblock Nodal auxiliary space preconditioning in h(curl) and h(div) spaces.
\newblock {\em {SIAM} Journal on Numerical Analysis}, 45(6):2483--2509, January
  2007.

\bibitem{Vassilevski2009}
T.V. Kolev and P.S. Vassilevski.
\newblock Parallel auxiliary space {AMG} for h(curl) problems.
\newblock {\em Journal of Computational Mathematics}, 27(5):604--623, June
  2009.

\bibitem{NECHAEV2008}
O.V. Nechaev, E.P. Shurina, and M.A. Botchev.
\newblock Multilevel iterative solvers for the edge finite element solution of
  the 3d maxwell equation.
\newblock {\em Computers {\&} Mathematics with Applications},
  55(10):2346--2362, 2008.
\newblock Advanced Numerical Algorithms for Large-Scale Computations.

\bibitem{Arnold2000}
D.N. Arnold, R.S. Falk, and R.~Winther.
\newblock Multigrid in h (div) and h (curl).
\newblock {\em Numerische Mathematik}, 85(2):197--217, April 2000.

\bibitem{NullSpace_Maxwell}
L.~Xue and D.~Jiao.
\newblock Method for analytically finding the nullspace of stiffness matrix for
  both zeroth-order and higher order curl-conforming vector bases in
  unstructured meshes.
\newblock {\em IEEE Transactions on Microwave Theory and Techniques},
  68(2):456--468, 2020.

\bibitem{Tikhonov1948}
A.N. Tikhonov and A.A. Samarskii.
\newblock On the radiation principle.
\newblock {\em Zh. Eksper. i Teoret. Fiz.}, 18(2):243--248, 1948.

\bibitem{Tikhonov1990}
A.N. Tikhonov and A.A. Samarskii.
\newblock {\em Equations of Mathematical Physics}.
\newblock Dover Publications, 1990.

\bibitem{Sommerfeld1912}
A.~Sommerfeld.
\newblock Die greensche funktion der schwingungslgleichung.
\newblock {\em Jahresbericht der Deutschen Mathematiker-Vereinigung},
  21:309--352, 1912.

\bibitem{Bristeau1998}
M.O. Bristeau, R.~Glowinski, and J.~P{\'{e}}riaux.
\newblock Controllability methods for the computation of time-periodic
  solutions: Application to scattering.
\newblock {\em Journal of Computational Physics}, 147(2):265--292, December
  1998.

\bibitem{Grote2019}
M.J. Grote and J.~H. Tang.
\newblock On controllability methods for the helmholtz equation.
\newblock {\em Journal of Computational and Applied Mathematics}, 358:306--326,
  October 2019.

\bibitem{ChaumontFrelet2022}
T.~Chaumont-Frelet, M.J. Grote, S.~Lanteri, and J.H. Tang.
\newblock A controllability method for maxwell{\textquotesingle}s equations.
\newblock {\em {SIAM} Journal on Scientific Computing}, 44(6):A3700--A3727,
  December 2022.

\bibitem{Appel2020}
D.~Appel\"{o}, F.~Garcia, and O.~Runborg.
\newblock {WaveHoltz}: Iterative solution of the helmholtz equation via the
  wave equation.
\newblock {\em {SIAM} Journal on Scientific Computing}, 42(4):A1950--A1983,
  January 2020.

\bibitem{Peng2022}
Z.~Peng and D.~Appelo.
\newblock {EM}-{WaveHoltz}: A flexible frequency-domain method built from
  time-domain solvers.
\newblock {\em {IEEE} Transactions on Antennas and Propagation},
  70(7):5659--5671, July 2022.

\bibitem{Courant1928}
R.~Courant, K.~Friedrichs, and H.~Lewy.
\newblock Über die partiellen differenzengleichungen der mathematischen
  physik.
\newblock {\em Mathematische Annalen}, 100:32--74, 1928.

\bibitem{Mikhailenko2008}
B.~G. Mikhailenko and A.~F. Mastryukov.
\newblock Numerical solution of maxwell's equations for anisotropic media using
  the laguerre transform.
\newblock {\em Russian Geology and Geophysics}, 49:621--627, 2008.

\bibitem{Mastryukov2013}
A.F. Mastryukov.
\newblock Solving an inverse problem for maxwell's equations numerically with
  laguerre functions.
\newblock {\em Numerical Analysis and Applications}, 6(4):279--288, October
  2013.

\bibitem{fatab2011}
A.~G. Fatyanov and A.~V. Terekhov.
\newblock High-performance modeling acoustic and elastic waves using the
  parallel dichotomy algorithm.
\newblock {\em J. Comp. Phys.}, 230(5):1992--2003, 2011.

\bibitem{Terekhov2015206}
A.~V. Terekhov.
\newblock Spectral-difference parallel algorithm for the seismic forward
  modeling in the presence of complex topography.
\newblock {\em Journal of Applied Geophysics}, 115(0):206--219, 2015.

\bibitem{Mikhailenko1999}
B.~G. Mikhailenko.
\newblock Spectral {Laguerre} method for the approximate solution of time
  dependent problems.
\newblock {\em Applied Mathematics Letters}, 12:105--110, 1999.

\bibitem{Terekhov2017}
A.~V. Terekhov.
\newblock The {Laguerre} finite difference one-way equation solver.
\newblock {\em Computer Physics Communications}, 214:71 -- 82, 2017.

\bibitem{Terekhov2018}
A.~V. Terekhov.
\newblock The stabilization of high-order multistep schemes for the {Laguerre}
  one-way wave equation solver.
\newblock {\em Journal of Computational Physics}, 368:115 -- 130, 2018.

\bibitem{NIST:DLMF}
{\it NIST Digital Library of Mathematical Functions}.
\newblock http://dlmf.nist.gov/, Release 1.1.2 of 2021-06-15.
\newblock F.~W.~J. Olver, A.~B. {Olde Daalhuis}, D.~W. Lozier, B.~I. Schneider,
  R.~F. Boisvert, C.~W. Clark, B.~R. Miller, B.~V. Saunders, H.~S. Cohl, and
  M.~A. McClain, eds.

\bibitem{Weeks1966}
William~T. Weeks.
\newblock Numerical inversion of {L}aplace transforms using {L}aguerre
  functions.
\newblock {\em J. ACM}, 13(3):419--429, July 1966.

\bibitem{Keilson1981}
J.~Keilson, W.~Nunn, and U.~Sumita.
\newblock The bilateral {Laguerre} transform.
\newblock {\em Applied Mathematics and Computation}, 8(2):137 -- 174, 1981.

\bibitem{Weber1980}
H.~Weber.
\newblock Numerical computation of the {F}ourier transform using laguerre
  functions and the fast {Fourier} transform.
\newblock {\em Numerische Mathematik}, 36(2):197--209, Jun 1980.

\bibitem{Integral_Transform}
L.~Debnath and D.~Bhatta.
\newblock {\em Integral Transforms and Their Applications, Second Edition}.
\newblock Taylor \& Francis, 2006.

\bibitem{Prabhu2018}
K.~M.~M. Prabhu.
\newblock {\em Window functions and their applications in signal processing}.
\newblock CRC Press, 2018.

\bibitem{Saad}
Y.~Saad.
\newblock {\em Iterative Methods for Sparse Linear Systems}.
\newblock SIAM, 2003.

\bibitem{Saad1986}
Y.~Saad and M.H. Schultz.
\newblock {GMRES}: A generalized minimal residual algorithm for solving
  nonsymmetric linear systems.
\newblock {\em {SIAM} Journal on Scientific and Statistical Computing},
  7(3):856--869, July 1986.

\end{thebibliography}
\end{document}